\documentclass[12pt]{article}%
\usepackage{amsmath}
\usepackage{amsfonts}
\usepackage{amssymb}
\usepackage{graphicx}%
\providecommand{\U}[1]{\protect\rule{.1in}{.1in}}
\newtheorem{theorem}{Theorem}

\newtheorem{remark}[theorem]{Remark}

\begin{document}

\title{Global bifurcation of vortex and dipole solutions in Bose-Einstein condensates}
\author{Andres Contreras \thanks{{\small Science Hall 224, New Mexico State
University, Department of Mathematical Sciences, USA. acontre@nmsu.edu}},
Carlos Garc\'{\i}a-Azpeitia \thanks{{\small Departamento de Matem\'{a}ticas,
Facultad de Ciencias, Universidad Nacional Aut\'{o}noma de M\'{e}xico, 04510
M\'{e}xico DF, M\'{e}xico. cgazpe@ciencias.unam.mx}}}
\maketitle

\begin{abstract}
The Gross-Pitaevskii equation for a Bose-Einstein condensate (BEC) with
symmetric harmonic trap is given in \eqref{Ec}. Periodic solutions of
\eqref{Ec} play an important role in the understanding of the long term
behavior of its solutions. In this note we prove the existence of several
global branches of solutions to \eqref{Ec} among which there are vortex
solutions and dipole solutions.

\end{abstract}

\section{Introduction}

The Gross-Pitaevskii equation for a Bose-Einstein condensate (BEC) with
symmetric harmonic trap is given by
\begin{equation}
-iu_{t}-\Delta u+(x^{2}+y^{2})u+\left\vert u\right\vert ^{2}u=0. \label{Ec}%
\end{equation}

Periodic solutions of \eqref{Ec} play an important role in the understanding
of the long term behavior of its solutions. In \cite{PeKe13}, symmetric and
asymmetric vortex solutions are obtained and their stability is established.
Solutions with two rotating vortices of opposite vorticity are constructed in
\cite{57}. In \cite{GKC} the authors prove the existence of periodic and
quasi-periodic trajectories of dipoles in anisotropic condensates. The
literature of the study of vortex dynamics in Bose-Einstein condensates is
vast, both on the mathematical and physical side; we refer the reader to
\cite{17,21,34,38,PeKe11,57} and the references therein for a more detailed account.

In this note we prove the existence of several global branches of solutions to
\eqref{Ec} among which there are vortex solutions and dipole solutions. Let
$X$ be the space of functions in $H^{2}(\mathbb{R}^{2};\mathbb{C)}$ for which
$\left\Vert u\right\Vert _{X}^{2}=\left\Vert u\right\Vert _{H^{2}}%
^{2}+\left\Vert r^{2}u\right\Vert _{L^{2}}^{2}$ is finite. Our main results are:

\begin{theorem}
\label{teo1} Let $m_{0}\geq1$ and $n_{0}$ be fixed non-negative integers. The
equation \eqref{Ec} has a global bifurcation in
\begin{equation}
Fix(\tilde{O}(2))=\{u\in X:u(r,\theta)=e^{im_{0}\theta}u(r)\text{ with
}u(r)\text{ real valued}\text{ }\}\text{.} \label{deffixteo1}%
\end{equation}
These are periodic solutions to \eqref{Ec} of the form
\[
e^{-i\omega t}e^{im_{0}\theta}u(r)\text{,}%
\]
starting from $\omega=2(m_{0}+2n_{0}+1)$, where $u(r)$ is a real-valued function.
\end{theorem}

It will be clear from the proof, and standard bifurcation theory, that for
small amplitudes $a$, we have the local expansion
\[
u(r)=av_{m_{0},n_{0}}(r)+\sum_{n\in\mathbb{N}}u_{m_{0},n}v_{m_{0},n}(r)\text{
and }u_{m_{0},n}(a)=O(a^{2})\text{,}%
\]
where the $v_{m,n}$'s are the eigenfunctions introduced in \eqref{def.vmn} and
the $u_{m,n}$'s are Fourier coefficients. Thus, the number $m_{0}$ is the
degree of the vortex at the origin and $n_{0}$ is the number of nodes of
$u(r)$ in $(0,\infty)$.

Our second theorem is concerned with the existence of multi-pole solutions.

\begin{theorem}
\label{teo2} Let $m_{0}\geq1$ and $n_{0}<m_{0}$ be two fixed non-negative
integers, then the equation \eqref{Ec} has a global bifurcation in
\begin{equation}
Fix(\mathbb{Z}_{2}\times\tilde{D}_{2m_{0}})=\{u\in X:u(r,\theta)=\bar
{u}(r,\theta)=u(r,-\theta)=-u(r,\theta+\pi/m_{0})\}\text{.} \label{deffixteo2}%
\end{equation}
These are periodic solutions to \eqref{Ec} of the form
\[
e^{-i\omega t}u(r,\theta)\text{,}%
\]
starting from $\omega=2(m_{0}+2n_{0}+1)$, where $u(r,\theta)$ is a real
function vanishing at the origin, enjoying the symmetries%
\[
u(r,\theta)=u(r,-\theta)=-u(r,\theta+\pi/m_{0})\text{.}%
\]

\end{theorem}

The requirement $n_{0}<m_{0}$ in Theorem \ref{teo2} is a non-resonance
condition. The solutions of the previous theorem for $(m_{0},n_{0})=(1,0)$
correspond to dipole solutions. This follows locally from the estimate
\[
u(r,\theta)=a(e^{im_{0}\theta}+e^{-im_{0}\theta})v_{m_{0},n_{0}}(r)+\sum
_{m\in\{m_{0},3m_{0},5m_{0},...\}}\sum_{n\in\mathbb{N}}u_{m,n}(e^{im\theta
}+e^{-im\theta})v_{m,n}(r)
\]
for small amplitude $a$ where $u_{m,n}=O(a^{2}).$

For $n_{0}=0$, since the function $v_{m_{0},0}$ is positive for $r\in
(0,\infty)$ and $u_{m,n}=O(a^{2})$, $u(r,\theta)$ is zero only when
$\theta=(k+1/2)\pi/m_{0}$. Moreover, as $u$ is real, then the lines
$\theta=(k+1/2)\pi/m_{0}$ correspond to zero density regions and the phase has
a discontinuous jump of $\pi$ at those lines.

A difficulty when trying to obtain the dipole solution, that is for
$(m_{0},n_{0})=(1,0)$, is the fact that to carry out a local inversion one has
to deal with a linearized operator with a repeated eigenvalue corresponding to
$(m_{0},n_{0})=(1,0)$ and $(-1,0)$. We overcome this by restricting our
problem to a natural space of symmetries which we identify below. In this
space our linearized operator only encounters a simple bifurcation which
yields the global existence result thanks to a topological degree argument;
see Theorem \ref{Theo}.

\section{Reduction to a bifurcation in a subspace of symmetries}

The group of symmetries of \eqref{Ec} is a three torus, corresponding to
rotations, phase and time invariances. The analysis of the group
representations leads to two kinds of isotropy groups, one corresponding to
vortex solutions and the other to dipole solutions. A fixed point argument on
restricted subspaces and Leray-Schauder degree yield the global existence of
these branches.

In \cite{PeKe13}, the authors study the case $(m_{0},n_{0})=(1,0)$ which is
bifurcation of a vortex of degree one. They also obtain second branch stemming
from this one and analyze its stability. The bifurcation from the case
$(m_{0},n_{0})=(0,0)$ is the ground state (see \cite{PeKe13}). The proof of
existence of dipole-like solutions was left open. In the present work we use
the symmetries of the problem, classifying the spaces of irreducible
representations, to obtain these as global branches provided a non-resonance
condition is satisfied.

\subsection{Setting the problem}

We rewrite \eqref{Ec} it in polar coordinates:
\[
-iu_{t}-\Delta_{(r,\theta)}u+r^{2}u+\left\vert u\right\vert ^{2}u=0\text{,}%
\]
where $\Delta_{(r,\theta)}=\partial_{r}^{2}+r^{-1}\partial_{r}+r^{-2}%
\partial_{\theta}^{2}$. Periodic solutions of the form $u(t,r,\theta
)=e^{-i\omega t}u(r,\theta)$ are zeros of the map
\[
f(u,\omega)=-\Delta_{(r,\theta)}u+(r^{2}-\omega+\left\vert u\right\vert
^{2})u.
\]

The eigenvalues and eigenfunctions of the linear Schr\"{o}dinger operator
\[
L=-\Delta_{(r,\theta)}+r^{2}:X\rightarrow L^{2}(\mathbb{R}^{2};\mathbb{C}),
\]
are found in Chapter 6, complement D, pg.727-737, on \cite{QM}. The operator
$L$ has eigenfunctions $v_{m,n}(r)e^{im\theta}$, which form and orthonormal
basis of $L^{2}(\mathbb{R}^{2};\mathbb{C}),$ and eigenvalues
\[
\lambda_{m,n}=2(\left\vert m\right\vert +2n+1)\text{ for }(m,n)\in
\mathbb{Z\times N}\text{,}%
\]
where $v_{m,n}(r)$ is a solution of $\left(  -\Delta_{m}+r^{2}\right)
v_{m,n}(r)=\lambda_{m,n}v_{m,n}(r)$, where%
\begin{equation}
-\Delta_{m}+r^{2}:=-(\partial_{r}^{2}+r^{-1}\partial_{r}-r^{-2}m^{2}%
)+r^{2}\text{,} \label{def.vmn}%
\end{equation}
with $v_{m,n}(0)=0$ for $m\neq0$. We have that
\[
u=\sum_{(m,n)\in\mathbb{Z\times N}}u_{m,n}v_{m,n}(r)e^{im\theta},\qquad
Lu=\sum_{(m,n)\in\mathbb{Z\times N}}\lambda_{m,n}u_{m,n}v_{m,n}(r)e^{im\theta
}\text{.}%
\]
Moreover, we know that $v_{m,n}(r)e^{im\theta}$ are orthogonal functions,
where $n$ is the number of nodes of $v_{m,n}(r)$ in $(0,\infty)$, see section
2.9 in \cite{Mo}.

\begin{remark}
Notice, that this is a slightly different orthonormal system that the one in
\cite{PeKe13}, which is more suited for anisotropic traps: $V(x,y)=\alpha
x^{2}+\beta y^{2}$ with $\alpha\neq\beta$.
\end{remark}

We have that the norm of $u$ in $L^{2}(\mathbb{R}^{2};\mathbb{C})$ is
$\left\Vert u\right\Vert _{L^{2}}^{2}=\sum_{(m,n)\in\mathbb{Z\times N}%
}\left\vert u_{m,n}\right\vert ^{2}$. Then, the inverse operator
$K=L^{-1}:L^{2}(\mathbb{R}^{2};\mathbb{C})\rightarrow X$ is continuous and
given by
\[
Ku=\sum\lambda_{m,n}^{-1}u_{m,n}v_{m,n}(r)e^{im\theta}.
\]
Moreover, the operator $K:X\rightarrow X$ is compact.

Observe that $H^{2}(\mathbb{R}^{2})$ is a Banach algebra and $H^{2}%
(\mathbb{R}^{2})\subset C^{0}(\mathbb{R}^{2})$, then $\left\Vert uv\right\Vert
_{X}\leq c\left\Vert v\right\Vert _{X}\left\Vert u\right\Vert _{X}$. We then
see that $g(u):=K(\left\vert u\right\vert ^{2}u)=\mathcal{O}(\left\Vert
u\right\Vert _{X}^{3})$ is a \textit{nonlinear compact map} such that
$g:X\rightarrow X$. Therefore, we obtain an equivalent formulation for the
bifurcation as zeros of the map%
\[
Kf(u,\omega)=u-\omega Ku+g(u):X\times\mathbb{R}\rightarrow X\text{.}%
\]
This formulation has the advantage that allows us to appeal to the global
Rabinowitz alternative \cite{Ra} (Theorem \ref{Theo} below).

\subsection{Equivariant bifurcation}

Let us define the action of the group generated by $(\psi,\varphi
)\in\mathbb{T}^{2}$, $\kappa\in\mathbb{Z}_{2}$ and $\bar{\kappa}\in
\mathbb{Z}_{2}$ in $L^{2}$ as%

\[
\rho(\psi,\varphi)u(r,\theta)=e^{i\varphi}u(r,\theta+\psi);\quad\rho
(\kappa)u(r,\theta)=u(r,-\theta);\quad\rho(\bar{\kappa})u(r,\theta)=\bar
{u}(r,\theta)\text{.}%
\]
Actually, the group generated by these actions is $\Gamma=O(2)\times O(2),$
and the map $Kf$ is $\Gamma$-equivariant.

Given a pair $(m_{0},n_{0})\in\mathbb{Z\times N},$ the operator $K$ has
multiple eigenvalues $\lambda_{m,n}^{-1}=\lambda_{m_{0},n_{0}}^{-1}$ for each
$(m,n)\in\mathbb{Z\times N}$ such that $\left\vert m\right\vert +2n=\left\vert
m_{0}\right\vert +2n_{0}$. To reduce the multiplicity of the eigenvalue
$\lambda_{m_{0},n_{0}}^{-1}$, we assume for the moment that there is a
subgroup $G$ of $\Gamma$ such that in the fixed point space,%
\[
Fix(G)=\{u\in X:\rho(g)u=u\text{ for }g\in G\}\text{,}%
\]
the linear map $K$ has only one eigenvalue $\lambda_{m_{0},n_{0}}^{-1}$. Then,
we can apply the following theorem using the fact that $Kf(u,\omega
):Fix(G)\times\mathbb{R}\rightarrow Fix(G)$ is well defined.

\begin{theorem}
\label{Theo}There is a global bifurcating branch $Kf(u(\omega),\omega)=0,$
starting from $\omega=\lambda_{m_{0},n_{0}}$ in the space $Fix(G)\times
\mathbb{R}$, this branch is a continuum that is unbounded or returns to a
different bifurcation point $(0,\omega_{1})$.
\end{theorem}

For a proof see the simplified approach due to Ize in Theorem 3.4.1 of
\cite{Ni2001}, or a complete exposition in \cite{IzVi03}.

We note that if $u\in H^{2}$ is a zero of $Kf$, then $u=K(\omega u-\left\vert
u\right\vert ^{2}u)\in H^{4}$. Using a bootstrapping argument we obtain that
the zeros of $Kf$ are solutions of the equation (\ref{Ec}) in $C^{\infty}.$
Next, we find the irreducible representations and the maximal isotropy groups.
The fixed point spaces of the maximal isotropy groups will have the property
that $K$ has a simple eigenvalue corresponding to $\lambda_{m_{0},n_{0}}%
^{-1}.$ Thus, Theorems \ref{teo1} and \ref{teo2} will follow from Theorem
\ref{Theo} applied to $G=\tilde{O}(2)$ and $G=\mathbb{Z}_{2}\times\tilde
{D}_{2m_{0}}$ respectively.

\subsection{Isotropy groups}

The action of the group on the components $u_{m,n}$ is given by $\rho
(\varphi,\psi)u_{m,n}=e^{i\varphi}e^{im\psi}u_{m,n}$ and $\rho(\kappa
)u_{m,n}=u_{-m,n}\text{ and }\rho(\bar{\kappa})u_{m,n}=\bar{u}_{-m,n}$. Then,
the irreducible representations are $(z_{1},z_{2})=(u_{m,n},u_{-m,n}%
)\in\mathbb{C}^{2}$, and the action of $\Gamma$ in the representation
$(z_{1},z_{2})$ is
\begin{align}
\rho(\varphi,\psi)(z_{1},z_{2})  &  =e^{i\varphi}(e^{im\psi}z_{1},e^{-im\psi
}z_{2});\label{Action}\\
\rho(\kappa)(z_{1},z_{2})  &  =(z_{2},z_{1});\nonumber\\
\rho(\bar{\kappa})(z_{1},z_{2})  &  =(\bar{z}_{2},\bar{z}_{1}).\nonumber
\end{align}

Actually, the irreducible representations $(u_{m_{0},n},u_{-m_{0},n}%
)\in\mathbb{C}^{2}$ are similar for all $n\in\mathbb{N}.$ The spaces of
similar irreducible representations are of infinite dimension. We analyze only
non-radial bifurcations, that is solutions bifurcating from $\omega
=\lambda_{m_{0},n_{0}}$ with $m_{0}\neq0 ;$ the radial bifurcation with
$m_{0}=0$ may be analyzed directly from the operator associated to the
spectral problem \eqref{def.vmn}.

Let us fix $m_{0}\geq1$ and $(z_{1},z_{2})=(u_{m_{0},n},u_{-m_{0},n})$. Then,
possibly after applying $\kappa$, we may assume $z_{1}\neq0$, unless
$(z_{1},z_{2})=(0,0)$. Moreover, using the action of $S^{1}$, the point
$(z_{1},z_{2})$ is in the orbit of $(a,re^{i\theta})$. It is known that there
are only two maximal isotropy groups, one corresponding to $(a,0)$ and the
other one to $(a,a)$, see for instance \cite{GoSc86}.

From (\ref{Action}), we have that the isotropy group of $(a,0)$ is generated
by $(\varphi,-\varphi/m_{0})$ and $\kappa\bar{\kappa}$, that is
\[
\tilde{O}(2)=\left\langle (\varphi,-\varphi/m_{0}),\kappa\bar{\kappa
}\right\rangle \text{.}%
\]
While the isotropy group of $(a,a)$ is generated by$\ (\pi,\pi/m_{0})$,
$\kappa$ and $\bar{\kappa}$, that is
\[
\mathbb{Z}_{2}\times\tilde{D}_{2m_{0}}=\left\langle \kappa,(\pi,\pi
/m_{0}),\bar{\kappa}\right\rangle \text{.}%
\]
These two groups are the only maximal isotropy groups of the representation
$(z_{1},z_{2})\in\mathbb{C}^{2}$, and the fixed point spaces have real
dimension one in $\mathbb{C}^{2}$.

\section{Vortex solutions: Proof of Theorem \ref{teo1}}

The functions fixed by the group $\tilde{O}(2)$ satisfy $u(r,\theta
)=e^{im_{0}\theta}u(r)$ from the element $(\varphi,-\varphi/m_{0})$, and
$u(r)=\bar{u}(r)$ from the element $\kappa\bar{\kappa}$. Thus, functions in
the space $Fix(\tilde{O}(2))$ are of the form
\[
u(r,\theta)=\sum_{n\in\mathbb{N}}u_{m_{0},n}e^{im_{0}\theta}v_{m_{0},n}(r)
\]
with $u_{m_{0},n}\in\mathbb{R}$.

Therefore, the map $Kf(u,\omega)$ has a simple eigenvalue $\lambda
_{m_{0},n_{0}}$ in the space $Fix(\tilde{O}(2))\times\mathbb{R}$ (see
\eqref{deffixteo1}). Therefore, from Theorem \ref{Theo}, there is a global
bifurcation in $Fix(\tilde{O}(2))\times\mathbb{R}$ starting at $\omega
=\lambda_{m_{0},n_{0}}.$

\section{Multi-pole-like solutions: Proof of Theorem \ref{teo2}}

The functions fixed by $\mathbb{Z}_{2}\times\tilde{D}_{2m_{0}}$ satisfy
$u(r,\theta)=\bar{u}(r,\theta)=u(r,-\theta)$. Therefore $u_{m,n}$ is real and
$u_{m,n}=u_{-m,n}$. Moreover, since $u(r,\theta)=-u(r,\theta+\pi/m_{0})$, then
$u_{m,n}=-e^{i\pi(m/m_{0})}u_{m,n}$. This relation gives $u_{m,n}=0$ unless
$e^{i\pi(m/m_{0})}=-1$ or $m/m_{0}$ is odd. Thus, functions in the space
$Fix(\mathbb{Z}_{2}\times\tilde{D}_{2m_{0}})$ are of the form
\[
u(r,\theta)=\sum_{m\in\{m_{0},3m_{0},5m_{0},...\}}%
{\displaystyle\sum\limits_{n\in\mathbb{N}}}
u_{m,n}(e^{im\theta}+e^{-im\theta})v_{m,n}(r)\text{,}%
\]
where $u_{m,n}$ is real and $m_{0}\geq1$. Therefore, the map $K$ has a simple
eigenvalue $\lambda_{m_{0},n_{0}}=2(m_{0}+2n_{0}+1)$ in $Fix(\mathbb{Z}%
_{2}\times\tilde{D}_{2m_{0}})$ if $\lambda_{lm_{0},n}\neq\lambda_{m_{0},n_{0}%
}$ for $n\in\mathbb{N}$ and $l=3,5,7..$.This condition is equivalent to
$lm_{0}+2n\neq m_{0}+2n_{0}$ or $2n_{0}-(l-1)m_{0}\neq2n$. Then, the
eigenvalue $\lambda_{m_{0},n_{0}}$ is simple if $2n_{0}<(l-1)m_{0}$ for
$l=3,5,...$, or $n_{0}<m_{0}$.

From Theorem (\ref{Theo}), the map $Kf(u,\omega)$ has a global bifurcation in
$Fix(\mathbb{Z}_{2}\times\tilde{D}_{2m_{0}})\times\mathbb{R}$ as $\omega$
crosses the value $\lambda_{m_{0},n_{0}}$ (see \textit{\eqref{deffixteo2}}).

\textbf{Acknowledgements.} C. Garcia learned about this problem from D.
Pelinovsky. A. Contreras is grateful to D.Pelinovsky for useful discussions.
Both authors thank the referees for their useful comments which greatly
improved the presentation of this note.

\end{document}